\documentstyle{amsppt}
\voffset-10mm
\magnification1200
\pagewidth{130mm}
\pageheight{204mm}
\hfuzz=2.5pt\rightskip=0pt plus1pt
\binoppenalty=10000\relpenalty=10000\relax
\TagsOnLeft
\loadbold
\nologo
\addto\tenpoint{\normalbaselineskip=1.2\normalbaselineskip\normalbaselines}
\addto\eightpoint{\normalbaselineskip=1.05\normalbaselineskip\normalbaselines}
\let\le\leqslant
\let\ge\geqslant
\let\[\lfloor
\let\]\rfloor
\let\ka=k
\redefine\d{\roman d}
\redefine\Re{\operatorname{Re}}
\define\sect#1#2#3{\parshape=3 4mm 117mm 9mm 112mm 9mm 112mm
\noindent\eightpoint{\bf#1}\enspace#2\dotfill\dotfill
\noindent\vskip-\normalbaselineskip\parshape=1 4mm 122mm {}
\hfill\eightpoint#3}
\topmatter
\title
Euler's constant, $q$-logarithms, and \\
formulas of Ramanujan and Gosper
\endtitle
\author
Jonathan Sondow {\rm(New York) and}
Wadim Zudilin\footnotemark"$^\ddag$"\ {\rm(Moscow)}
\endauthor
\date
26 March 2003
\enddate
\address
\hbox to70mm{\vbox{\hsize=70mm%
\leftline{209 West 97th Street}
\leftline{New York, NY 10025 USA}
\leftline{{\it URL\/}: \tt http://home.earthlink.net/\~{}jsondow/}
}}
\endaddress
\email
{\tt jsondow\@alumni.princeton.edu}
\endemail
\address
\hbox to70mm{\vbox{\hsize=70mm%
\leftline{Moscow Lomonosov State University}
\leftline{Department of Mechanics and Mathematics}
\leftline{Vorobiovy Gory, GSP-2, Moscow 119992 RUSSIA}
\leftline{{\it URL\/}: \tt http://wain.mi.ras.ru/index.html}
}}
\endaddress
\email
{\tt wadim\@ips.ras.ru}
\endemail
\abstract
The aim of the paper is to relate computational and
arithmetic questions about Euler's constant $\gamma$ with properties
of the values of the $q$-logarithm function,
with natural choice of~$q$. By these means,
we generalize a classical formula for $\gamma$
due to Ramanujan, together with Vacca's and
Gosper's series for~$\gamma$, as well as deduce
irrationality criteria and tests and new asymptotic
formulas for computing Euler's constant. The main tools
are Euler-type integrals and hypergeometric series.
\endabstract
\endtopmatter

\noindent
\footnotetext"$^\ddag$"{The work of the second author
is supported by an Alexander von Humboldt research fellowship.}
\footnote""{{\it Key words and phrases}.\enspace
Euler's constant, $q$-logarithm, hypergeometric series,
Euler's transform, irrationality.}
\footnote""{2000 {\it Mathematics Subject Classification}.\enspace
Primary 11Y60; Secondary 11J72, 33C20, 33D15.}
\rightheadtext{Euler's constant, $q$-logarithms and Ramanujan}
\leftheadtext{J.~Sondow and W.~Zudilin}
\document

We recall the definition of Euler's constant:
$$
\split
\gamma
&=\lim_{n\to\infty}\biggl(\sum_{\nu=1}^n\frac1\nu-\log n\biggr)
\\
&=0.57721566490153286060651209008240243104215933593992\dots
\endsplit
$$
and the fact that the (expected) irrationality of~$\gamma$
has not yet been proved.

Recently, the first author gave~\cite{9}, \cite{10} a
construction of $\Bbb Z$-linear forms involving~$\gamma$
and logarithms; namely, he proved that
$$
d_{2n}I_n
\in\Bbb Z+\Bbb Z\gamma
+\Bbb Z\log(n+1)+\Bbb Z\log(n+2)+\dots+\Bbb Z\log(2n),
$$
where
$$
\align
I_n
&=\iint\limits_{[0,1]^2}
\frac{x^n(1-x)^ny^n(1-y)^n}{(1-xy)|\log xy|}\,\d x\,\d y
\\
&=\sum_{\nu=n+1}^\infty\int_\nu^\infty
\biggl(\frac{n!}{t(t+1)\dotsb(t+n)}\biggr)^2\d t
\endalign
$$
and $d_n$ denotes the least common multiple of the numbers $1,2,\dots,n$.
This type of approximation allows deducing
several (ir)rationality criteria for Euler's constant~\cite{9}.
In a sense, what follows is inspired by considerations in~\cite{9}.

A special function that we will require is the $q$-logarithm
$$
\ln_q(1+z)=\sum_{\nu=1}^\infty\frac{(-1)^{\nu-1}z^\nu}{q^\nu-1},
\qquad |z|<q,
\tag1
$$
for $q$~real with $|q|>1$ (in fact, our main choices
will be $q=2$ in Sections~2--4 and $q=3$ in Section~8,
but we let $q>1$ be any integer in Section~9).
The function~\thetag{1} is a $q$-extension of the ordinary logarithm
function
$$
\log(1+z)=\sum_{\nu=1}^\infty\frac{(-1)^{\nu-1}z^\nu}\nu,
\qquad |z|\le1, \quad z\ne-1,
$$
since
$$
\lim\Sb q\to1\\q>1\endSb\bigl((q-1)\ln_q(1+z)\bigr)=\log(1+z),
\qquad |z|<1.
$$
We also mention that the irrationality of the values of
the $q$-logarithm~\thetag{1} for
$q\in\Bbb Z\setminus\{0,\allowmathbreak\pm1\}$
and $z\in\Bbb Q$, $z\le1$, is known~\cite{3}.

\head
Contents
\endhead
\eightpoint

\sect{\phantom01.}{Gosper's acceleration of Vacca's series}{2}

\sect{\phantom02.}{A $q$-logarithm approach}{4}

\sect{\phantom03.}{Irrationality tests for Euler's constant}{8}

\sect{\phantom04.}{Computing Euler's constant}{10}

\sect{\phantom05.}{Ramanujan's base $q$ integral for Euler's constant}{11}

\sect{\phantom06.}{A base $q$ generalization of Gosper's Series}{12}

\sect{\phantom07.}{Base $q$ integrals}{13}

\sect{\phantom08.}{Base $3$ irrationality testing}{15}

\sect{\phantom09.}{Base $q$ asymptotic formulas for Euler's constant}{18}

\sect{10.}{Concluding remarks}{22}

\sect{\phantom{00.}}{References}{23}

\tenpoint

\head
1. Gosper's acceleration of Vacca's series
\endhead

The following expression for~$\gamma$ is known
as Vacca's series~\cite{11}:
$$
\gamma
=\sum_{n=1}^\infty(-1)^n\frac{\[\log_2n\]}n,
\tag2
$$
where $\log_qx=(\log x)/(\log q)$ is the ordinary logarithm base~$q$.
Gosper \cite{6} has transformed this series
into a more rapidly converging series of positive rationals.
His method (see Section~6) is to use partial summation
to write \thetag{2} as the double series
$$
\gamma
=\sum_{\nu=1}^\infty\sum_{\ka=0}^\infty\frac{(-1)^\ka}{2^\nu+\ka},
\tag3
$$
apply Euler's transformation to each inner sum on~$\ka$,
obtaining the double series~\thetag{7} below,
then triangularize it. In this way,
he arrives at a version of series~\thetag{4}.
Here we give a different proof of~\thetag{4}.
To begin, we reverse the order of summation in~\thetag{3} and
write the resulting double series as the integral~\thetag{5},
which is a form of an integral that Bromwich
attributes to Catalan (\cite{4}, p.~526).

\proclaim{Theorem 1 \rm(Gosper's Series)}
We have the series for Euler's constant
$$
\gamma
=\frac12+\sum_{\nu=2}^\infty\frac1{2^{\nu+1}}
\sum_{\ka=1}^{\nu-1}\binom{2^{\nu-\ka}+\ka}\ka^{-1}.
\tag4
$$
\endproclaim

\demo{Proof}
After expanding $1/(1+x)$ in a geometric series, termwise integration
shows that the formula
$$
\gamma
=\int_0^1\frac1{1+x}\sum_{\nu=1}^\infty x^{2^\nu-1}\,\d x
\tag5
$$
is equivalent to~\thetag{3} with the order of summation reversed,
which is easily justified by grouping terms in pairs.
Integrating by parts $K+1$ times gives
$$
\gamma
=\sum_{\ka=0}^K\sum_{\nu=1}^\infty
\frac1{2^{\nu+\ka+1}\binom{2^\nu+\ka}\ka}
+\int_0^1\frac{(K+1)!}{(1+x)^{K+2}}
\sum_{\nu=1}^\infty\frac{x^{2^\nu+K}}{2^\nu(2^\nu+1)\dotsb(2^\nu+K)}\,\d x.
\tag6
$$
Since the integral is less than
$$
\int_0^1\frac1{(1+x)^{K+2}}
\sum_{\nu=1}^\infty\frac1{2^\nu}\,\d x
=\frac1{K+1}\biggl(1-\frac1{2^{K+1}}\biggr),
$$
as $K$ tends to infinity \thetag{6} becomes
$$
\gamma
=\sum_{\ka=0}^\infty\sum_{\nu=1}^\infty
\frac1{2^{\nu+\ka+1}\binom{2^\nu+\ka}\ka}.
\tag7
$$
The inner series with $\ka=0$ sums to~$1/2$,
and in the series with $\ka>0$ we may
collect terms (triangularize) as in~\thetag{4},
completing the proof. \qed
\enddemo

\remark{Remark}
We can regard the series on the right
of~\thetag{2},~\thetag{4},~\thetag{7}
as base $q=2$ series for Euler's constant.
Gosper indicated~\cite{6} the following extensions
of~\thetag{2},~\thetag{7} to base $q=3$:
$$
\gamma
=2\sum_{n=1}^\infty\frac{\[\log_3n\]\cdot\cos\frac{2\pi n}3}n
=2\sum_{\nu=1}^\infty\sum_{k=0}^\infty
\frac{\cos\frac{\pi(k-1)}6}{3^{\nu+(k+1)/2}\binom{3^\nu+k}k}
\tag8
$$
(see Sections 5 and 6 for proofs and generalizations),
which give faster convergence than formulas~\thetag{2},~\thetag{7},
respectively.
\endremark

\head
2. A $q$-logarithm approach
\endhead

\proclaim{Lemma 1}
Let $q$ be real, $q>1$.
Then for all integers $n\ge0$ and~$\ka$, with $0<\ka<q^{n+1}$,
the following identity holds:
$$
\sum_{\nu=n+1}^\infty\frac1{q^\nu+\ka}
=\frac1\ka\,\ln_q\biggl(1+\frac\ka{q^n}\biggr).
\tag9
$$
\endproclaim

\remark{Remark}
For $\ka=0$, the series in \thetag{9} equals
$$
\sum_{\nu=n+1}^\infty\frac1{q^\nu}
=\frac{q^{-(n+1)}}{1-q^{-1}}
=\frac1{q^n(q-1)}.
$$
\endremark

\demo{Proof}
We have
$$
\allowdisplaybreaks
\align
\sum_{\nu=n+1}^\infty\frac1{q^\nu+\ka}
&=\sum_{\nu=n+1}^\infty\frac1{q^\nu}\,\frac1{1+\ka\cdot q^{-\nu}}
=\sum_{\nu=n+1}^\infty\frac1{q^\nu}
\sum_{\mu=1}^\infty\frac{(-\ka)^{\mu-1}}{q^{\nu(\mu-1)}}
\\
&=\sum_{\mu=1}^\infty(-\ka)^{\mu-1}
\sum_{\nu=n+1}^\infty\frac1{q^{\nu\mu}}
=\sum_{\mu=1}^\infty(-\ka)^{\mu-1}\frac{q^{-(n+1)\mu}}{1-q^{-\mu}}
\\
&=\frac1\ka\sum_{\mu=1}^\infty
\frac{(-1)^{\mu-1}(\ka/q^n)^\mu}{q^\mu-1}
=\frac1\ka\,\ln_q\biggl(1+\frac\ka{q^n}\biggr)
\endalign
$$
as required. \qed
\enddemo

\remark{Remark}
P.~Sebah~\cite{8} has pointed out that in order to compute
the $q$-logarithm, the following formula is much more
efficient than~\thetag{1} or~\thetag{9}:
$$
\ln_q(1+z)
=z\sum_{\nu=1}^N\frac1{q^\nu+z}+r_N(z),
\qquad\text{where}\quad
r_N(z)=\sum_{\nu=1}^\infty
\frac{(-1)^{\nu-1}z^\nu}{q^{N\nu}(q^\nu-1)},
$$
with $N$ being any positive integer.
The proof of the formula is similar to that of Lemma~1.
\endremark

Let $m$ be a non-negative integer.
Using the rational function
$$
R_m(t)=\frac{m!}{t(t+1)\dotsb(t+m)}
=\sum_{j=0}^m\frac{(-1)^j\binom mj}{t+j},
\tag10
$$
for each $\nu\ge0$ define the convergent series
$$
S_{\nu,m}
=\sum_{t=0}^\infty(-1)^tR_m(2^\nu+t).
\tag11
$$

\proclaim{Lemma 2}
We have
$$
S_{\nu,m}
=2^m\sum_{\ka=0}^\infty\frac{(-1)^\ka}{2^\nu+\ka}
-\sum_{j=1}^m\binom mj\sum_{\ka=0}^{j-1}
\frac{(-1)^\ka}{2^\nu+\ka}.
\tag12
$$
\endproclaim

(We use the standard convention that any sum $\sum_{j=a}^b$ is
zero if $a>b$.)

\demo{Proof}
Indeed,
$$
\align
S_{\nu,m}
&=\sum_{j=0}^m\sum_{t=0}^\infty\frac{(-1)^{t+j}\binom mj}{2^\nu+t+j}
=\sum_{j=0}^m\binom mj\sum_{t=0}^\infty\frac{(-1)^{t+j}}{2^\nu+t+j}
\\
&=\sum_{j=0}^m\binom mj\sum_{\ka=j}^\infty\frac{(-1)^\ka}{2^\nu+\ka}
=\sum_{j=0}^m\binom mj
\biggl(\sum_{\ka=0}^\infty-\sum_{\ka=0}^{j-1}\biggr)
\frac{(-1)^\ka}{2^\nu+\ka},
\endalign
$$
and \thetag{12} follows.
\qed
\enddemo

\proclaim{Lemma 3}
For the series $S_{\nu,m}$, we have the representations
$$
S_{\nu,m}
=(-1)^{2^\nu}\sum_{t=2^\nu}^\infty(-1)^tR_m(t)
=\int_0^1\frac{(1-x)^m}{1+x}x^{2^\nu-1}\,\d x.
\tag13
$$
\endproclaim

\demo{Proof}
Replace $t$ by $t-2^\nu$ in~\thetag{11} to get the first
equality in~\thetag{13}. For the second,
use the Binomial Theorem to write
$$
R_m(t)=\int_0^1(1-x)^mx^{t-1}\,\d x.
$$
Substitute this in the series for $S_{\nu,m}$, interchange summation
and integration, and sum the series
$$
\sum_{t=2^\nu}^\infty(-1)^{t-2^\nu}x^{t-1}
=\frac{x^{2^\nu-1}}{1+x},
\tag14
$$
arriving at the second equality in~\thetag{13}.
To justify termwise integration, replace the geometric
series~\thetag{14} with a finite sum plus remainder
$r_N(x)=(-1)^Nx^{N+2^\nu-1}/\allowmathbreak(1+x)$, and note that
$$
\biggl|\int_0^1r_N(x)\,\d x\biggr|<\frac1{N+2^\nu}\to0
\qquad\text{as $N\to\infty$}.
\qed
$$
\enddemo

Our final object is the double series
$$
I_{n,m}=\sum_{\nu=n+1}^\infty S_{\nu,m},
\tag15
$$
where $n\ge0$. By Lemma~2 and formula~\thetag{3},
for the quantity~\thetag{15} we have
$$
\split
I_{n,m}
&=2^m\biggl(\sum_{\nu=1}^\infty-\sum_{\nu=1}^n\biggr)
\sum_{\ka=0}^\infty\frac{(-1)^\ka}{2^\nu+\ka}
-\sum_{j=1}^m\binom mj\sum_{\ka=0}^{j-1}
\sum_{\nu=n+1}^\infty\frac{(-1)^\ka}{2^\nu+\ka}
\\
&=2^m\gamma
-2^m\sum_{\nu=1}^n
\biggl(\sum_{l=1}^\infty-\sum_{l=1}^{2^\nu-1}\biggr)\frac{(-1)^l}l
-\sum_{j=1}^m\binom mj\sum_{\nu=n+1}^\infty\frac1{2^\nu}
\\ &\qquad
-\sum_{j=1}^m\binom mj\sum_{\ka=1}^{j-1}
\sum_{\nu=n+1}^\infty\frac{(-1)^\ka}{2^\nu+\ka}.
\endsplit
$$
For simplicity, assume now that $m\le2^{n+1}$. Then using Lemma~1,
$$
\split
I_{n,m}
&=2^m\cdot(\gamma+n\log2)
-2^m\sum_{\nu=1}^n\sum_{l=1}^{2^\nu-1}\frac{(-1)^{l-1}}l
-\frac{2^m-1}{2^n}
\\ &\qquad
-\sum_{j=2}^m\binom mj\sum_{\ka=1}^{j-1}\frac{(-1)^\ka}\ka\,
\ln_2\biggl(1+\frac\ka{2^n}\biggr).
\endsplit
\tag16
$$
We can summarize the result of these computations
in the following way.

\proclaim{Lemma 4}
If $m\le2^{n+1}$, then
$$
I_{n,m}=2^m\gamma+L_{n,m}-A_{n,m},
\tag17
$$
where
$$
L_{n,m}
=2^mn\log2
+\sum_{j=2}^m\binom mj\sum_{\ka=1}^{j-1}\frac{(-1)^{\ka-1}}\ka\,
\ln_2\biggl(1+\frac\ka{2^n}\biggr)
$$
and
$$
A_{n,m}
=\frac{2^m-1}{2^n}
+2^m\sum_{\nu=1}^n\sum_{l=1}^{2^\nu-1}\frac{(-1)^{l-1}}l
$$
satisfy
$$
d_mL_{n,m}\in\Bbb Z\log2
+\Bbb Z\ln_2\biggl(1+\frac1{2^n}\biggr)
+\Bbb Z\ln_2\biggl(1+\frac2{2^n}\biggr)+\dots
+\Bbb Z\ln_2\biggl(1+\frac{m-1}{2^n}\biggr)
$$
and $d_{2^n}A_{n,m}\in\Bbb Z$.
\endproclaim

\proclaim{Lemma 5}
For the double series $I_{n,m}$, we have the integral
$$
I_{n,m}
=\int_0^1\frac{(1-x)^m}{1+x}
\sum_{\nu=n+1}^\infty x^{2^\nu-1}\,\d x,
\tag18
$$
which generalizes Catalan's integral~\thetag{5} for~$\gamma$.
\endproclaim

\demo{Proof}
Integrating termwise (since the series in the integrand is dominated
by a geometric series, the justification in the proof of Lemma~3
applies here) and using Lemma~3, the formula follows
from the definition of~$I_{n,m}$.
The case $I_{0,0}=\gamma$ (see~\thetag{16})
is indeed the integral~\thetag{5}.
\qed
\enddemo

\proclaim{Lemma 6}
For $m$ a positive integer and $r\ge1$ a real number, we have the bounds
$$
\int_0^1x^{rm-1}(1-x)^m\,\d x
>\frac1{rm}\biggl(\frac{r^r}{(r+1)^{r+1}}\biggr)^m
\tag19
$$
and, for $0<x<1$,
$$
x^{rm-1}(1-x)^m
<4\cdot\biggl(\frac{r^r}{(r+1)^{r+1}}\biggr)^m.
\tag20
$$
\endproclaim

\demo{Proof}
Euler's beta integral
$$
\int_0^1x^{\alpha-1}(1-x)^{\beta-1}\,\d x
=\frac{\Gamma(\alpha)\,\Gamma(\beta)}{\Gamma(\alpha+\beta)},
\qquad \Re\alpha>0, \quad \Re\beta>0,
$$
gives
$$
\int_0^1x^{rm-1}(1-x)^m\,\d x
=\frac{\Gamma(rm)\,m!}{\Gamma(rm+m+1)}
=\frac1{rm}\prod_{j=1}^m\biggl(1+\frac{rm}j\biggr)^{-1};
$$
denote the product by $\Pi(m,r)$. Using $r\ge1$ to get
$$
\frac{j+r(m+1)}{j+rm}
<1+\frac1m
\le\biggl(1+\frac1r\biggr)^{r/m}
\qquad\text{for}\quad j=1,2,\dots,m,
$$
we obtain
$$
\frac{\Pi(m+1,r)}{\Pi(m,r)}
=\frac1{r+1}\prod_{j=1}^m\frac{j+rm}{j+r(m+1)}
>\frac1{r+1}\biggl(1+\frac1r\biggr)^{-r}
=\frac{r^r}{(r+1)^{r+1}},
$$
and by induction on $m$ it follows that
$$
\Pi(m,r)>\biggl(\frac{r^r}{(r+1)^{r+1}}\biggr)^m
$$
since the inequality holds for $m=1$.
This proves~\thetag{19}.

If $0<x<1$, then
$$
\split
x^{rm-1}(1-x)^m
&<\bigl(x^r(1-x)\bigr)^{m-1/r}
\le\biggl(\frac{r^r}{(r+1)^{r+1}}\biggr)^{m-1/r}
\\
&=\frac{(r+1)^{1+1/r}}r
\cdot\biggl(\frac{r^r}{(r+1)^{r+1}}\biggr)^m
\endsplit
$$
and \thetag{20} follows, using $r\ge1$.
\qed
\enddemo

\proclaim{Lemma 7}
If $n\ge0$ is an integer and $r$
is a real number satisfying $1\le r\le 2^{n+1}$,
then we have the bounds
$$
\frac1{2r(m+1)}\biggl(\frac{r^r}{(r+1)^{r+1}}\biggr)^{m+1}
<I_{n,m}
<6\cdot\biggl(\frac{r^r}{(r+1)^{r+1}}\biggr)^m
\qquad\text{for}\quad
m=\biggl\[\frac{2^{n+1}}r\biggr\].
\tag21
$$
In particular, for $r\ge1$ fixed and
$m=\[2^{n+1}/r\]\ge1$, the quantity~$I_{n,m}$
is positive and decreases exponentially at the rate
$$
\lim_{n\to\infty}
\frac{\log I_{n,\[2^{n+1}/r\]}}{\[2^{n+1}/r\]}
=r\log r-(r+1)\log(r+1).
\tag22
$$
\endproclaim

\demo{Proof}
Put $m=\[2^{n+1}/r\]$ in Lemma~5. Then
$m+1>2^{n+1}/r$, hence
$$
I_{n,m}
>\int_0^1\frac{(1-x)^m}{1+x}x^{2^{n+1}-1}\,\d x
>\frac12\int_0^1x^{r(m+1)-1}(1-x)^{m+1}\,\d x,
$$
and the lower bound in~\thetag{21} follows from~\thetag{19}
(with $m+1$ in place of~$m$).

Using $m\le2^{n+1}/r$ and relations \thetag{20} and \thetag{5},
we have
$$
\split
I_{n,m}
&\le\int_0^1\frac{(1-x)^m}{1+x}x^{rm-1}
\sum_{\nu=n+1}^\infty x^{2^\nu-2^{n+1}}\,\d x
\\
&<\int_0^1x^{rm-1}(1-x)^m\frac1{1+x}
\biggl(1+\sum_{\nu=1}^\infty x^{2^\nu-1}\biggr)\,\d x
\\
&<4\biggl(\frac{r^r}{(r+1)^{r+1}}\biggr)^m(\log2+\gamma),
\endsplit
$$
and the upper bound in~\thetag{21} follows from
the inequality $4(\log2+\gamma)<6$.
Since \thetag{21} implies \thetag{22}, we are done.
\qed
\enddemo

\remark{Remark \rom1}
If $2^{n+1}/r$ is an integer, then we may replace $m+1$ by $m$ in the
lower bound in~\thetag{21}.
\endremark

\remark{Remark \rom2}
As $r\to\infty$, the right-hand side of~\thetag{22} equals
$-\log r+O(1)$, hence its absolute value tends to infinity with~$r$.
\endremark

\head
3. Irrationality tests for Euler's constant
\endhead

\proclaim{Theorem 2 \rm(Rationality Criterion for $\gamma$)}
The fractional part of $d_{2^n}L_{n,m}$ equals
$d_{2^n}I_{n,m}$ for all $n$ and $m$ with
$$
\text{$n$~sufficiently large and $2^{n-1}\le m\le2^{n+1}$},
\tag23
$$
if and only if Euler's constant is a rational number.
\endproclaim

\demo{Proof}
With $2^{n-1}\le m\le2^{n+1}$ set $r=2^{n+1}/m$, so that $1\le r\le4$.
By Lemma~7, we have $I_{n,m}<6\rho(r)^m=6\rho(r)^{2^n\cdot2/r}$,
where $\rho(r)=r^r/(r+1)^{r+1}$. Straightforward verification shows
that $\rho(r)^{2/r}$ is increasing for $r>0$, and satisfies
$$
0<\rho(r)^{2/r}<\frac1e
\qquad\text{for}\quad 0<r<r_0=5.6213305349\dots,
$$
where $\rho(r_0)^{2/r_0}=1/e$.
Therefore
$I_{n,m}<6\bigl(\rho(4)^{1/2}\bigr)^{2^n}$
and $\rho(4)^{1/2}<1/e$.
Since the Prime Number Theorem implies that for any $\varepsilon>0$
we have $d_N=O\bigl((e^{1+\varepsilon})^N\bigr)$ as $N\to\infty$,
we conclude that $0<d_{2^n}I_{n,m}<1$ for $n$ and $m$
satisfying~\thetag{23}.

Now multiply \thetag{17} by $d_{2^n}$ and write the result as
$$
d_{2^n}L_{n,m}=d_{2^n}I_{n,m}+d_{2^n}(A_{n,m}-2^m\gamma).
\tag24
$$
Since $d_{2^n}A_{n,m}\in\Bbb Z$ by Lemma~4, the theorem follows.
\qed
\enddemo

\remark{Remark}
The argument shows that the Criterion holds with \thetag{23}
replaced by the more general conditions $n\ge n(r_1)$ and
$\lceil2^{n+1}/r_1\rceil\le m\le2^{n+1}$, where $r_1$~is any
number between~$1$ and~$r_0$.
\endremark

\proclaim{Theorem 3 \rm(Irrationality Test for $\gamma$)}
If the fractional part of $d_{2^n}L_{n,2^n}$ satisfies
$$
\bigl\{d_{2^n}L_{n,2^n}\bigr\}
>6\cdot\biggl(\frac{128}{729}\biggr)^{2^{n-1}}
\tag25
$$
infinitely often, then Euler's constant is irrational.
In fact, the inequality for a given~$n>0$ implies that if $\gamma$
is a rational number, then its denominator does not divide~$2^{2^n}d_{2^n}$.
\endproclaim

\demo{Proof}
It follows from \cite{7}, Theorem~13, that $d_{2N}<8^N$
for all integers $N\ge1$ (see \cite{9}, Lemma~3).
Combined with the upper bound in Lemma~7 for $r=2$, this gives
$$
0<d_{2^n}I_{n,2^n}
<8^{2^{n-1}}\cdot6\cdot\biggl(\frac4{27}\biggr)^{2\cdot2^{n-1}}
=6\cdot\biggl(\frac{128}{729}\biggr)^{2^{n-1}}
$$
for $n\ge1$. The first part of the theorem follows by the
Rationality Criterion for~$\gamma$, and the second part by
formula~\thetag{24} with $m=2^n$.
\qed
\enddemo

\remark{Remark \rom1}
Condition~\thetag{25} in Theorem~3 can be replaced
by the simpler condition
$$
\bigl\{d_{2^n}L_{n,2^n}\bigr\}
>5^{-2^{n-1}}
$$
if $n>1$. For $n>4$, this follows from the inequalities
$128/729=0.17558299\ldots<1/5$ and
$6\cdot(128/729)^{14}<5^{-14}$.
For $n=2,3,4$, it follows from the proof of Theorem~3
by replacing the bound $d_{2m}<8^m$ with
the exact values $d_4=12$, $d_8=840$, $d_{16}=720720$,
for which $d_{2^n}\cdot6\cdot(4/27)^{2^n}<5^{-2^{n-1}}$.
\endremark

\remark{Remark \rom2}
Using the Rationality Criterion for~$\gamma$, we may
generalize the Irrationality Test as follows:
{\it $\gamma$~is irrational if there exist infinitely
many integers $n$ and $m$ such that $2^{n-1}\le m\le2^{n+1}$ and}
$$
\bigl\{d_{2^n}L_{n,m}\bigr\}
>8^{2^{n-1}}\cdot6\cdot\biggl(\frac{r^r}{(r+1)^{r+1}}\biggr)^m,
\qquad\text{where}\quad
r=\frac{2^{n+1}}m.
\tag26
$$
For example, the case $m=2^n$ is Theorem~3, and the case $m=2^{n+1}$
can be stated: {\it a sufficient condition for irrationality of~$\gamma$
is that}
$$
\bigl\{d_{2^n}L_{n,2^{n+1}}\bigr\}
>6\cdot2^{-5\cdot2^{n-1}}
\qquad\text{\it infinitely often}.
$$
For any $\varepsilon>0$, we may further refine the Test by
replacing $8^{2^{n-1}}\cdot6$ with $e^{(1+\varepsilon)2^n}$
in~\thetag{26}.
\endremark

\example{Example}
P.~Sebah~\cite{8} has calculated that
$$
\bigl\{d_{2^{12}}L_{12,2^{12}}\bigr\}
=0.178346164\dotsc.
$$
Since the last number exceeds $5^{-2^{11}}$,
{\it if $\gamma$ is a rational number $a/b$, then
$b$~is not a divisor of $2^{4096}\cdot d_{4096}$\rom;
in particular, $|b|\ge4099$\/} (since the numbers
$4097=17\cdot241$ and $4098=2\cdot3\cdot683$ divide $d_{4096}$).
\endexample

\remark{Remark}
Similar irrationality tests for $\gamma$, but using ordinary
logarithms instead of $2$@-logarithms, are proved in~\cite{9},
where an example shows that no divisor of
\linebreak
$\binom{20000}{10000}\cdot d_{20000}$ can be a denominator of~$\gamma$.
The present example yields additional information because
of the high power of~$2$.
\endremark

\head
4. Computing Euler's constant
\endhead

\proclaim{Theorem 4}
The following asymptotic formulas as $n\to\infty$ are valid:
$$
\align
\gamma&=\frac{A_{n,2^{n+1}}-L_{n,2^{n+1}}}{2^{2^{n+1}}}
+O(2^{-6\cdot2^n}),
\tag27
\\
\gamma&=\frac{A_{n,2^n}-L_{n,2^n}}{2^{2^n}}
+O\Biggl(\biggl(\frac2{27}\biggr)^{2^n}\Biggr),
\tag28
\\
\gamma&=\frac{A_{n,2^{n-1}}-L_{n,2^{n-1}}}{2^{2^{n-1}}}
+O\Biggl(\biggl(\frac{128}{3125}\biggr)^{2^{n-1}}\Biggr),
\tag29
\\
\gamma&=\frac12A_{n,1}-n\log2+O(2^{-n}).
\tag30
\endalign
$$
\endproclaim

\demo{Proof}
For \thetag{27}, put $m=2^{n+1}$ in~\thetag{17} and use the upper bound
in~\thetag{21} with $r=1$.
For~\thetag{28}, \thetag{29}, take $m=2^n,2^{n-1}$ and
$r=2,4$, respectively.
Finally, let $m=1$ in Lemma~4 to get $I_{n,1}=2\gamma+2n\log2-A_{n,1}$.
From~\thetag{21} with $r=2^{n+1}$, we obtain
$$
0<I_{n,1}<6\cdot\frac{(2^{n+1})^{2^{n+1}}}{(2^{n+1}+1)^{2^{n+1}+1}}
=\frac6{2^{n+1}+1}\biggl(1+\frac1{2^{n+1}}\biggr)^{-2^{n+1}}
<2^{-n+1}
$$
for $n\ge0$, proving more than \thetag{30}.
\qed
\enddemo

\remark{Remark}
Using \thetag{28} and his calculation of $L_{n,2^n}$ for $n=12$,
P.~Sebah~\cite{8} has computed $\gamma$ correct to 4631~decimal places,
which is the accuracy $(2/27)^{2^{12}}$ predicted.
Formula~\thetag{27} gives almost 60\% more digits of~$\gamma$,
but requires computing about twice as many $2$-logarithms.
Similarly, \thetag{29} yields nearly 39\% fewer digits than~\thetag{28},
but involves only half as many $2$-logarithms. Letting $m$ continue
down to~$1$, we reach formula~\thetag{30}, which doesn't involve
any $2$-logarithms.
\endremark

\head
5. Ramanujan's base $q$ integral for Euler's constant
\endhead

For integer $q\ge2$, Ramanujan
(see~\cite{2} for proofs and references)
gave the formula for Euler's constant
$$
\gamma=\int_0^1\biggl(\frac q{1-x^q}-\frac1{1-x}\biggr)
\sum_{\nu=1}^\infty x^{q^\nu-1}\,\d x,
\tag31
$$
which reduces to Catalan's integral~\thetag{5} when $q=2$.
Equation~\thetag{31} implies the following base $q$
generalization of Vacca's series~\thetag{2}
(see \cite{2}, Theorem~2.6):
$$
\gamma
=\sum_{n=1}^\infty\frac{\sigma_n\[\log_qn\]}n,
\tag32
$$
where
$$
\sigma_n=\sigma_{n,q}
=\cases
q-1 &\text{if $q\mid n$}, \\
-1 &\text{if $q\nmid n$}.
\endcases
\tag33
$$
The case $q=3$ gives the first equality in~\thetag{8}.

Let $\epsilon=\epsilon_q$ be a fixed primitive $q$-th root of unity
(for instance, $\epsilon=e^{2\pi i/q}$).
Using
$$
\sigma_n
=-1+\sum_{l=0}^{q-1}\epsilon^{nl}
=\sum_{l=1}^{q-1}\epsilon^{nl},
$$
we can represent formula~\thetag{32} as
$$
\gamma
=\sum_{l=1}^{q-1}\sum_{n=1}^\infty\frac{\[\log_qn\]\epsilon^{nl}}n,
\tag34
$$
which can be regarded as another base $q$
generalization of Vacca's series~\thetag{2}.

\head
6. A base $q$ generalization of Gosper's Series
\endhead

Following Gosper's proof~\cite{6} of~\thetag{7},
we write \thetag{34} as
$$
\gamma
=\sum_{l=1}^{q-1}\sum_{\nu=1}^\infty
\sum_{k=0}^\infty\frac{\epsilon^{kl}}{q^\nu+k}
=\sum_{l=1}^{q-1}\sum_{\nu=1}^\infty\frac1{q^\nu}
\cdot{}_2F_1\biggl(\matrix 1, \; q^\nu \\ q^\nu+1 \endmatrix
\biggm|\epsilon^l\biggr)
$$
and apply Euler's transform
(see, e.g., \cite{1}, Section~2.4, formula~(1))
$$
{}_2F_1\biggl(\matrix\format\c&\,\c\\
\alpha, & \beta \\ & \gamma \endmatrix
\biggm|z\biggr)
=\frac1{(1-z)^\alpha}
\cdot{}_2F_1\biggl(\matrix\format\c&\,\c\\
\alpha, & \gamma-\beta \\ & \gamma \endmatrix
\biggm|\frac{-z}{1-z}\biggr),
$$
yielding
$$
\split
\gamma
&=\sum_{l=1}^{q-1}
\sum_{\nu=1}^\infty\frac1{q^\nu}\cdot\frac1{1-\epsilon^l}
\cdot{}_2F_1\biggl(\matrix 1, \; 1 \\ q^\nu+1 \endmatrix
\biggm|\frac{-\epsilon^l}{1-\epsilon^l}\biggr)
\\
&=\sum_{l=1}^{q-1}
\sum_{\nu=1}^\infty\frac1{q^\nu}\sum_{k=0}^\infty
{\binom{q^\nu+k}k}^{-1}\frac{(-1)^k\epsilon^{kl}}
{(1-\epsilon^l)^{k+1}}.
\endsplit
$$
Substituting $x=1$ into the expansion
$$
\prod_{j=1}^{q-1}(1-\epsilon^jx)
=\frac{1-x^q}{1-x}=1+x+x^2+\dots+x^{q-1}
$$
leads to
$$
\frac1{1-\epsilon^l}
=\frac1q\cdot\prod\Sb j=1\\j\ne l\endSb^{q-1}(1-\epsilon^j)
\tag35
$$
(where the empty product in the case $q=2$ equals~$1$).
Therefore, we arrive at the following result.

\proclaim{Theorem 5}
We have the base $q\ge2$ accelerated series
$$
\gamma
=\sum_{\nu=1}^\infty\sum_{k=0}^\infty
\frac{(-1)^k\chi_q(k)}{q^{\nu+k+1}\binom{q^\nu+k}k},
\tag36
$$
where
$$
\chi_q(k)
=\sum_{l=1}^{q-1}\epsilon^{kl}
\prod\Sb j=1\\j\ne l\endSb^{q-1}(1-\epsilon^j)^{k+1},
$$
a symmetric polynomial in the roots of the polynomial $1+x+x^2+\dots+x^q$,
is an integer-valued function independent of the choice
of primitive $q$-th root of unity~$\epsilon$.
\endproclaim

\remark{Remark}
From \thetag{35} we have
$$
\epsilon^{kl}
\prod\Sb j=1\\j\ne l\endSb^{q-1}(1-\epsilon^j)^{k+1}
=\epsilon^{kl}\frac{q^{k+1}}{(1-\epsilon^l)^{k+1}},
$$
and it follows that
$$
|\chi_q(k)|
\le(q-1)\cdot\frac{q^{k+1}}{|1-e^{2\pi i/q}|^{k+1}}
=(q-1)\biggl(\frac q{2\sin\frac\pi q}\biggr)^{k+1}.
$$
Thus we have the relation
$$
\limsup_{k\to\infty}|\chi_q(k)|^{1/k}
\le\frac q{2\sin\frac\pi q},
$$
which shows how the convergence of series~\thetag{34}
has been accelerated in~\thetag{36}.
\endremark

\example{Examples}
Series \thetag{36} is a generalization of
Gosper's series~\thetag{7}, which is the case $q=2$.
When $q=3$, taking $\epsilon=e^{2\pi i/3}$ we obtain
$$
\split
\chi_3(k)
&=e^{2\pi ki/3}(1-e^{-2\pi i/3})^{k+1}
+e^{-2\pi ki/3}(1-e^{2\pi i/3})^{k+1}
\\
&=e^{2\pi ki/3}\cdot3^{(k+1)/2}e^{\pi(k+1)i/6}
+e^{-2\pi ki/3}\cdot3^{(k+1)/2}e^{-\pi(k+1)i/6}
\\
&=3^{(k+1)/2}\cdot2\cos\frac{\pi(k-1)}6,
\endsplit
$$
which proves the second equality in~\thetag{8}.
When $q=4$, we take $\epsilon=e^{\pi i/2}=i$ and obtain
$$
\split
\chi_4(k)
&=i^k2^{k+1}(1+i)^{k+1}
+(-1)^k(1-i)^{k+1}(1+i)^{k+1}
+(-i)^k2^{k+1}(1-i)^{k+1}
\\
&=2^{3(k+1)/2+1}\cdot\cos\frac{\pi(3k+1)}4
+(-1)^k2^{k+1},
\endsplit
$$
which is an integer-valued function of~$k$.
\endexample

\head
7. Base $q$ integrals
\endhead

There are several ways to generalize integrals~\thetag{18}
for $I_{n,m}$ and \thetag{31} for~$\gamma$ simultaneously.
One way is to define
$$
I_{n,m,q}
=\int_0^1\biggl(\frac q{1-x^q}-\frac1{1-x}\biggr)(1-x)^m
\sum_{\nu=n+1}^\infty x^{q^\nu-1}\,\d x.
\tag37
$$
Then $I_{n,m,2}=I_{n,m}$ and $I_{0,0,q}=\gamma$,
and the analog of Lemma~7 with $r=q$ and $m=q^n$,
$$
\lim_{n\to\infty}\frac{\log I_{n,q^n,q}}{q^n}
=q\log q-(q+1)\log(q+1),
\tag38
$$
can be easily derived by the methods of Section~2
(cf\. Lemma~9 below).
We may expand integral~\thetag{37} as in
the proof of Lemma~4. Namely, with $R_m(t)$ defined in~\thetag{10},
and $\sigma_t$ in~\thetag{33}, we have
$$
\frac q{1-x^q}-\frac1{1-x}
=\frac{(q-1)-x-x^2-\dotsb-x^{q-1}}{1-x^q}
=\sum_{t=0}^\infty\sigma_tx^t,
$$
hence
$$
\split
I_{n,m,q}
&=\sum_{t=0}^\infty\sigma_t\sum_{\nu=n+1}^\infty
\int_0^1x^{q^\nu+t-1}(1-x)^m\,\d x
=\sum_{\nu=n+1}^\infty\sum_{t=0}^\infty\sigma_tR_t(q^\nu+t)
\\
&=\sum_{j=0}^m(-1)^j\binom mj
\sum_{\nu=n+1}^\infty\sum_{k=j}^\infty
\frac{\sigma_{k-j}}{q^\nu+k}.
\endsplit
$$
Unfortunately, the last representation shows that $I_{n,m,q}$
involves not only logarithms and $q$-logarithms
but also ``generalized Euler constants''
$$
\gamma_{j,q}=\sum_{\nu=1}^\infty\sum_{k=j}^\infty
\frac{\sigma_{k-j}}{q^\nu+k},
\qquad j=0,1,\dots,q-1;
\tag39
$$
the equality $\gamma_{0,q}=\gamma$ follows from~\thetag{32}
just as \thetag{3} follows from Vacca's series~\thetag{2}
(the case $q=2$).
We may avoid the appearance of generalized Euler constants
$\gamma_{j,q}\ne\gamma$
in the particular case $q=3$ (details are in Section~8),
but not for $q>3$ (see Section~10).

However, there is a way to generalize integral~\thetag{18}
for $I_{n,m}$ that solves this problem, namely, by defining
$$
\split
I_{n,m,q}'
&=\int_0^1\biggl(\frac q{1-x^q}-\frac1{1-x}\biggr)
\biggl(q-\frac{1-x^q}{1-x}\biggr)^m
\sum_{\nu=n+1}^\infty x^{q^\nu-1}\,\d x
\\
&=\int_0^1\frac{(q-1-x-x^2-\dotsb-x^{q-1})^{m+1}}{1-x^q}
\sum_{\nu=n+1}^\infty x^{q^\nu-1}\,\d x,
\endsplit
\tag40
$$
with the properties $I_{n,m,2}'=I_{n,m}$ and $I_{0,0,q}'=\gamma$,
as for~\thetag{37}. For integral~$I_{n,m,q}'$ in~\thetag{40},
we may generalize the results of Sections~2 and~4,
but the poor asymptotics of the integral make these
generalizations useless for irrationality testing.
For example, the choice $m=(q^n-1)/(q-1)$ gives the inclusions
$$
\split
d_{q^n}I_{n,m,q}'
\in\Bbb Z+\Bbb Z\gamma+\Bbb Z\log q
&
+\Bbb Z\ln_q\biggl(1+\frac1{q^n}\biggr)
+\Bbb Z\ln_q\biggl(1+\frac2{q^n}\biggr)+\dotsb
\\ &\qquad
+\Bbb Z\ln_q\biggl(1+\frac{q^n-2}{q^n}\biggr)
\endsplit
$$
(see Lemma~10 below) and the asymptotics
$$
\lim_{n\to\infty}\frac{\log I_{n,(q^n-1)/(q-1),q}}{q^n}>-1
$$
(see Lemma 13);
this makes it impossible to use $I_{n,m,q}'$ to give irrationality
tests for $q>2$ similar to those of Theorems~2 and 3 for $q=2$.
Instead we use $I_{n,m,q}'$
to extend Theorem~4 to a base $q>2$ asymptotic formula for~$\gamma$;
details are in Section~9.

\head
8. Base $3$ irrationality testing
\endhead

Taking $q=3$ in~\thetag{31} gives one of Ramanujan's
formulas for Euler's constant
$$
\gamma=\int_0^1\frac{2+x}{1+x+x^2}
\sum_{\nu=1}^\infty x^{3^\nu-1}\,\d x
\tag41
$$
(see \cite{2}, Corollary 2.3 and the equation following (2.7)),
and taking $q=3$ in~\thetag{37} gives the base~$3$ integral
$$
I_{n,m,3}
=\int_0^1\frac{2+x}{1+x+x^2}(1-x)^m
\sum_{\nu=n+1}^\infty x^{3^\nu-1}\,\d x.
\tag42
$$
It turns out that, for $m$ a multiple of~$6$, we can extend
the results of Sections 2--4 to base $q=3$, as follows.

Since
$$
\frac{(1-x)^6}{1+x+x^2}
=28-34x+21x^2-7x^3+x^4-\frac{27}{1+x+x^2},
$$
by induction we have
$$
\frac{(1-x)^{6k}}{1+x+x^2}
=Q_k(x)+\frac{(-3)^{3k}}{1+x+x^2},
\tag43
$$
where $Q_k(x)$ is the polynomial defined by the recursion
$$
\gathered
Q_0(x)=0,
\\
Q_{k+1}(x)=(1-x)^6Q_k(x)+(-3)^{3k}(28-34x+21x^2-7x^3+x^4)
\quad\text{for}\; k\ge0.
\endgathered
\tag44
$$
From \thetag{42}, \thetag{43}, \thetag{41}, for $m=6k$ we obtain
$$
\split
I_{n,m,3}
&=\int_0^1(2+x)\biggl(Q_k(x)+\frac{(-3)^{3k}}{1+x+x^2}\biggr)
\biggl(\sum_{\nu=1}^\infty-\sum_{\nu=1}^n\biggr)x^{3^\nu-1}\,\d x
\\
&=(-3)^{3k}\biggl(\gamma-\int_0^1\frac{2+x}{1+x+x^2}
\sum_{\nu=1}^nx^{3^\nu-1}\,\d x\biggr)
\\ &\qquad
+\int_0^1(2+x)Q_k(x)\sum_{\nu=n+1}^\infty x^{3^\nu-1}\,\d x.
\endsplit
\tag45
$$

By division, we find that
$$
\frac{(2+x)x^{3^\nu-1}}{1+x+x^2}
=1+x-\sum_{\mu=1}^{3^{\nu-1}-1}
(2x^{3\mu-1}-x^{3\mu}-x^{3\mu+1})-\frac{1+2x}{1+x+x^2}
$$
for $\nu\ge1$, where the sum is zero if $\nu=1$. Therefore,
$$
\int_0^1\frac{2+x}{1+x+x^2}\sum_{\nu=1}^nx^{3^\nu-1}
=n\biggl(\frac32-\log3\biggr)
-\sum_{\nu=2}^n\sum_{\mu=1}^{3^{\nu-1}-1}
\biggl(\frac2{3\mu}-\frac1{3\mu+1}-\frac1{3\mu+2}\biggr)
\tag46
$$
for $n\ge1$, where the double sum vanishes if $n=1$.

To evaluate the last integral in~\thetag{45}, assume that
$m=6k<3^{n+1}$ and let
$a_{0,k},\dots,\allowmathbreak a_{6k-1,k}\in\Bbb Z$
denote the coefficients in the polynomial
$$
\sum_{j=0}^{6k-1}a_{j,k}x^j=(2+x)Q_k(x).
\tag47
$$
Noting that $a_{0,k}=(-3)^{3k}56k$ by~\thetag{44}, and using Lemma~1,
we have
$$
\split
\int_0^1(2+x)Q_k(x)\sum_{\nu=n+1}^\infty x^{3^\nu-1}\,\d x
&=\sum_{j=0}^{6k-1}a_{j,k}\sum_{\nu=n+1}^\infty\frac1{3^\nu+j}
\\
&=(-1)^k\frac{28k}{3^{n-3k}}
+\sum_{j=1}^{6k-1}\frac{a_{j,k}}j\ln_3\biggl(1+\frac j{3^n}\biggr).
\endsplit
\tag48
$$

Formulas \thetag{45}, \thetag{46}, \thetag{48} imply
the following base $3$ version of Lemma~4.

\proclaim{Lemma 8}
If $m<3^{n+1}$ is a multiple of~$6$, say $m=6k$, then
$$
I_{n,m,3}=(-3)^{3k}\gamma+L_{n,m,3}-A_{n,m,3},
\tag49
$$
where
$$
L_{n,m,3}
=(-3)^{3k}n\log3
+\sum_{j=1}^{m-1}\frac{a_{j,k}}j\ln_3\biggl(1+\frac j{3^n}\biggr)
$$
and
$$
A_{n,m,3}
=(-3)^{3k}\Biggl(\frac{3n}2-\frac{28k}{3^n}
-\sum_{\nu=2}^n\sum_{\mu=1}^{3^{\nu-1}-1}
\biggl(\frac2{3\mu}-\frac1{3\mu+1}-\frac1{3\mu+2}\biggr)\Biggr)
$$
satisfy
$$
d_mL_{n,m,3}\in\Bbb Z\log3
+\Bbb Z\ln_3\biggl(1+\frac1{3^n}\biggr)
+\Bbb Z\ln_3\biggl(1+\frac2{3^n}\biggr)+\dots
+\Bbb Z\ln_3\biggl(1+\frac{m-1}{3^n}\biggr)
$$
and $d_{3^n}A_{n,m,3}\in\Bbb Z$, and the $a_{j,k}$ are integers
determined by~\thetag{47} and~\thetag{44}.
\endproclaim

\proclaim{Lemma 9}
We have the bounds
$$
\frac1{40m+90}\biggl(\frac{3^3}{4^4}\biggr)^m
<I_{n,m,3}<3\cdot\biggl(\frac{3^3}{4^4}\biggr)^m
\qquad\text{for}\quad m=3^n-3\ge0
\tag50
$$
and the limit
$$
\lim_{n\to\infty}\frac{\log I_{n,3^n-3,3}}{3^n}
=3\log3-4\log4
=-2.24934057\dotsc.
\tag51
$$
\endproclaim

\demo{Proof}
Following the proof of Lemma~7, with $m=3^n-3$ we have
$$
\split
I_{n,m,3}
&>\int_0^1\frac{2+x}{1+x+x^2}(1-x)^mx^{3^{n+1}-1}\,\d x
>\int_0^1x^{3m+8}(1-x)^m\,\d x
\\
&=\frac{(3m+8)!\,m!}{(4m+9)!}
>\frac{(3/4)^8}{4m+9}\cdot{\binom{4m}m}^{-1}
>\frac1{40m+90}\biggl(\frac{3^3}{4^4}\biggr)^m
\endsplit
$$
and
$$
\split
I_{n,m,3}
&=\int_0^1\frac{2+x}{1+x+x^2}(1-x)^mx^{3m+8}
\sum_{\nu=n+1}^\infty x^{3^\nu-3^{n+1}}\,\d x
\\
&<\int_0^1x^{3m+8}(1-x)^m\frac{2+x}{1+x+x^2}
\biggl(1+\sum_{\nu=1}^\infty x^{3^\nu-1}\biggr)\,\d x
\\
&<\biggl(\frac{3^3}{4^4}\biggr)^m(2+\gamma)
<3\cdot\biggl(\frac{3^3}{4^4}\biggr)^m.
\endsplit
$$
This proves the lemma.
\qed
\enddemo

\proclaim{Theorem 6}
We have the asymptotic formula for Euler's constant
$$
\gamma=(-1)^{n-1}\frac{A_{n,3^n-3,3}-L_{n,3^n-3,3}}{3^{(3^n-3)/2}}
+O\Biggl(\biggl(\frac{3^{5/2}}{64}\biggr)^{3^n}\Biggr)
\qquad\text{as}\; n\to\infty.
$$
\endproclaim

\demo{Proof}
Let $m=6k=3^n-3$ for $n>0$. Then $(-1)^{3k}=(-1)^{n-1}$ and
the formula follows from \thetag{49} and~\thetag{51}.
\qed
\enddemo

\proclaim{Theorem 7}
Euler's constant is rational if and only if
the fractional part of
\linebreak
$d_{3^n}L_{n,3^n-3,3}$
equals $d_{3^n}I_{n,3^n-3,3}$ for all large~$n$.
\endproclaim

\demo{Proof}
Multiply \thetag{49} by $d_{3^n}$ and write the result as
$$
d_{3^n}L_{n,m,3}=d_{3^n}I_{n,m,3}+d_{3^n}(A_{n,m,3}-(-3)^{3k}\gamma).
$$
By the Prime Number Theorem and Lemma~9,
for $m=3^n-3$ and any $\varepsilon>0$ we have
$d_{3^n}I_{n,m,3}=O\bigl((3^3e^{1+\varepsilon}/4^4)^m\bigr)$ as $m\to\infty$.
Thus $0<d_{3^n}I_{n,3^n-3,3}<1$ for all large~$n$.
Since $d_{3^n}A_{n,3^n-3,3}\in\Bbb Z$,
the theorem follows.
\qed
\enddemo

\proclaim{Theorem 8}
If the inequality
$$
\bigl\{d_{3^n}L_{n,3^n-3,3}\bigr\}
>3\cdot\biggl(\frac34\biggr)^{3^{n+1}}
$$
holds infinitely often, then $\gamma$ is irrational.
In fact, the inequality for a given~$n>0$ implies that
no divisor of $3^{3^n}d_{3^n}$ can be a denominator of~$\gamma$.
\endproclaim

\demo{Proof}
The inequalities $d_{3N}\le d_{4N}<8^{2N}$ and the upper
bound in~\thetag{50} imply that
$$
d_{3^n}I_{n,3^n-3,3}
<8^{2\cdot3^{n-1}}\cdot3\cdot\biggl(\frac{3^3}{4^4}\biggr)^{3^n}
=3\cdot\biggl(\frac34\biggr)^{3^{n+1}}
$$
for $n\ge1$. Using Lemma~8,
this proves the second part of the theorem, which implies the first.
\qed
\enddemo

\head
9. Base $q$ asymptotic formulas for Euler's constant
\endhead

Denoting by $F_q(x)$, $G_q(x)$ the polynomials
$$
\gathered
F_q(x)=1+x+x^2+\dots+x^{q-1}=\frac{1-x^q}{1-x},
\\
G_q(x)=(q-1)+(q-2)x+\dots+x^{q-2}
=\frac{q-F_q(x)}{1-x},
\endgathered
\tag52
$$
we may write the integral~\thetag{40} as
$$
I_{n,m,q}'=\int_0^1\frac{G_q(x)}{F_q(x)}\bigl(q-F_q(x)\bigr)^m
\biggl(\sum_{\nu=1}^\infty-\sum_{\nu=1}^n\biggr)x^{q^\nu-1}\,\d x.
\tag53
$$
Expanding the binomial and using Ramanujan's formula $I_{0,0,q}=\gamma$,
we obtain
$$
\split
I_{n,m,q}'
&=q^m\biggl(\gamma-\int_0^1\frac{G_q(x)}{F_q(x)}
\sum_{\nu=1}^nx^{q^\nu-1}\,\d x\biggr)
\\ &\qquad
+\int_0^1G_q(x)\sum_{j=1}^m(-1)^j\binom mjq^{m-j}F_q(x)^{j-1}
\sum_{\nu=n+1}^\infty x^{q^\nu-1}\,\d x.
\endsplit
\tag54
$$

In the first integrand, we find by dividing that
$$
\frac{G_q(x)x^{q^\nu-1}}{F_q(x)}
=F_{q-1}(x)+\sum_{\mu=1}^{q^{\nu-1}-1}(F_q(x)-q)x^{q\mu-1}
-\frac{F_q'(x)}{F_q(x)}
$$
for $\nu\ge1$, the summation being zero if $\nu=1$. Therefore,
$$
\int_0^1\frac{G_q(x)}{F_q(x)}\sum_{\nu=1}^nx^{q^\nu-1}\,\d x
=n\biggl(\sum_{\mu=1}^{q-1}\frac1\mu-\log q\biggr)
-\sum_{\nu=2}^n\sum_{\mu=1}^{q^{\nu-1}-1}
\biggl(\frac{q-1}{q\mu}-\sum_{\lambda=1}^{q-1}\frac1{q\mu+\lambda}\biggr)
\tag55
$$
for $n\ge1$, where the double summation vanishes if $n=1$.

In the second integrand in~\thetag{54}, let
$b_{0,m,q},\dots,b_{m(q-1)-1,m,q}\in\Bbb Z$ denote the coefficients
of the polynomial
$$
\sum_{l=0}^{m(q-1)-1}b_{l,m,q}x^l
=G_q(x)\sum_{j=1}^m(-1)^j\binom mjq^{m-j}F_q(x)^{j-1}.
\tag56
$$
Note that $b_{0,m,q}=(q-1)\bigl((q-1)^m-q^m\bigr)$ by~\thetag{52}.
If $m\le q^{n+1}/(q-1)$, then by Lemma~1 the last integral in~\thetag{54}
equals
$$
\split
&
\int_0^1\sum_{l=0}^{m(q-1)-1}b_{l,m,q}x^l
\sum_{\nu=n+1}^\infty x^{q^\nu-1}\,\d x
=\sum_{l=0}^{m(q-1)-1}b_{l,m,q}\sum_{\nu=n+1}^\infty\frac1{q^\nu+l}
\\ &\qquad
=\frac{(q-1)^m-q^m}{q^n}
+\sum_{l=1}^{m(q-1)-1}\frac{b_{l,m,q}}l\,
\ln_q\biggl(1+\frac l{q^n}\biggr).
\endsplit
\tag57
$$

Equations \thetag{54}--\thetag{57} imply
the following generalization of Lemma~4 to base~$q$.

\proclaim{Lemma 10}
If $m\le q^{n+1}/(q-1)$, then
$$
I_{n,m,q}'=q^m\gamma+L_{n,m,q}'-A_{n,m,q}',
\tag58
$$
where
$$
L_{n,m,q}'
=q^mn\log q+\sum_{l=1}^{m(q-1)-1}\frac{b_{l,m,q}}l\,
\ln_q\biggl(1+\frac l{q^n}\biggr)
$$
and
$$
A_{n,m,q}'
=\frac{q^m-(q-1)^m}{q^n}+q^m\Biggl(\sum_{\mu=1}^{q-1}\frac n\mu
-\sum_{\nu=2}^n\sum_{\mu=1}^{q^{\nu-1}-1}\biggl(\frac{q-1}{q\mu}
-\sum_{\lambda=1}^{q-1}\frac1{q\mu+\lambda}\biggr)\Biggr)
$$
satisfy
$$
\split
d_{m(q-1)}L_{n,m,q}'
&\in\Bbb Z\log q
+\Bbb Z\ln_q\biggl(1+\frac1{q^n}\biggr)
+\Bbb Z\ln_q\biggl(1+\frac2{q^n}\biggr)+\dotsb
\\ &\qquad
+\Bbb Z\ln_q\biggl(1+\frac{m(q-1)-1}{q^n}\biggr)
\endsplit
$$
and $d_{q^n}A_{n,m,q}'\in\Bbb Z$, and the $b_{l,m,q}$ are integers
determined by \thetag{56} and~\thetag{52}.
\endproclaim

To optimize the asymptotics of $I_{n,m,q}'/q^m$
in the range of Lemma~10, and thus obtain the smallest error
term in our formula for~$\gamma$, we choose $m=(q^n-1)/(q-1)$
in the remainder of this section.
Let $f_q(x)$ denote the polynomial
$$
f_q(x)=\bigl(q-F_q(x)\bigr)x^{q(q-1)}
=(q-1)x^{q(q-1)}-\sum_{\nu=1}^{q-1}x^{q(q-1)+\nu}.
$$

\proclaim{Lemma 11}
If $q\ge2$, there exists $x_q$ between $0$ and~$1$ such that
$$
f_q(x_q)=\max_{0\le x\le1}f_q(x).
\tag59
$$
Moreover, we have the asymptotic formula
$$
\int_0^1f_q(x)^mx^{q-1}\,\d x=f_q(x_q)^{m(1+o(1))}
\qquad\text{as $m\to\infty$}
\tag60
$$
and the bounds
$$
\frac1{(q+1)e}<f_q(x_q)<\frac1{2e}.
\tag61
$$
\endproclaim

\demo{Proof}
Since
$$
\frac{f_q'(x)}{x^{q(q-1)-1}}
=q(q-1)^2-\sum_{\nu=1}^{q-1}(q(q-1)+\nu)x^\nu
\tag62
$$
is decreasing for $x\ge0$, positive at $x=0$ and
negative at $x=1$, it follows that $f_q'(x)$ vanishes
at exactly one point $x=x_q$ in the open interval $(0,1)$.
Then \thetag{59} holds with this choice of~$x_q$, since
$f_q(0)=f_q(1)=0$ and $f_q(x)>0$ for $0<x<1$.
Differentiating \thetag{62}, we substitute $f_q'(x_q)=0$ and
find that $f_q''(x_q)<0$; by Laplace's method
(see, e.g., \cite{5}, Section~4.2) we obtain \thetag{60}.

Since
$$
q-F_q(x)=(q-1)\biggl(1-\frac{x+x^2+\dots+x^{q-1}}{q-1}\biggr),
$$
the inequalities
$$
(x\cdot x^2\dotsb x^{q-1})^{1/(q-1)}
\le\frac{x+x^2+\dots+x^{q-1}}{q-1}
\le x
\qquad\text{for $0\le x\le1$}
$$
yield
$$
(q-1)(1-x)x^{q(q-1)}\le f_q(x)\le(q-1)(1-x^{q/2})x^{q(q-1)}.
$$
Using the relations
$$
\frac1{(r+1)e}
<\max_{0\le y\le1}\bigl((1-y)y^r\bigr)=\frac{r^r}{(r+1)^{r+1}}
<\frac1{re}
\qquad\text{for $r>0$},
$$
we set $y=x$, $r=q(q-1)$ and deduce the lower bound in~\thetag{61},
then set $y=x^{q/2}$, $r=2(q-1)$ and obtain the upper bound.
\qed
\enddemo

The next result with $q=2$ is comparable to Lemma~7 with $m=2^n-1$.

\proclaim{Lemma 12}
For $q\ge2$, we have the bounds
$$
0<I_{n,m,q}'<\frac q{(2e)^m}
\qquad\text{with}\quad m=\frac{q^n-1}{q-1}
\tag63
$$
and
$$
-1-\log(q+1)
<\lim_{n\to\infty}
\frac{\log I_{n,(q^n-1)/(q-1),q}'}{(q^n-1)/(q-1)}
<-1-\log2.
\tag64
$$
\endproclaim

\demo{Proof}
Let $m=(q^n-1)/(q-1)$, so that $x^{q^{n+1}-1}=x^{(mq+1)(q-1)}$.
Using~\thetag{53},~\thetag{59} and Ramanujan's formula
$I_{0,0,q}'=\gamma$, we have
$$
\align
0<I_{n,m,q}'
&=\int_0^1\frac{G_q(x)}{F_q(x)}\bigl(q-F_q(x)\bigr)^m
x^{(mq+1)(q-1)}\sum_{\nu=n+1}^\infty x^{q^\nu-q^{n+1}}\,\d x
\\
&<\int_0^1f_q(x)^m\frac{G_q(x)}{F_q(x)}
\biggl(1+\sum_{\nu=1}^\infty x^{q^\nu-1}\biggr)\,\d x
\tag65
\\
&\le f_q(x_q)^m(c_q+\gamma),
\qquad\text{where}\quad
c_q=\int_0^1\frac{G_q(x)}{F_q(x)}\,\d x,
\endalign
$$
and \thetag{63} follows using $c_q+\gamma<(q-1)+\gamma<q$
and the upper bound in~\thetag{61}.

We also have
$$
I_{n,m,q}'
>\int_0^1\frac{G_q(x)}{F_q(x)}\bigl(q-F_q(x)\bigr)^m
x^{(mq+1)(q-1)}\,\d x
>\frac{q-1}2\int_0^1f_q(x)^mx^{q-1}\,\d x.
\tag66
$$
From \thetag{65}, \thetag{66}, \thetag{60} we deduce that
$$
\lim_{n\to\infty}
\frac{\log I_{n,(q^n-1)/(q-1),q}'}{(q^n-1)/(q-1)}
=\log f_q(x_q),
\tag67
$$
and \thetag{61} gives \thetag{64}.
\qed
\enddemo

\remark{Remark}
The exact value of the constant $c_q$ is
$-\psi(1/q)-\log q-\gamma$, as another formula of Ramanujan
shows (see \cite{2}, Proposition~2.1).
\endremark

\proclaim{Theorem 9 \rm(Base $q$ Asymptotic Formula for $\gamma$)}
If $q\ge2$ and $m=(q^n-1)/(q-1)$, then
$$
\gamma=\frac{A_{n,m,q}'-L_{n,m,q}'}{q^m}+\delta_{m,q}
\qquad\text{with}\quad
0<\delta_{m,q}<\frac q{(2eq)^m}.
$$
\endproclaim

\demo{Proof}
Set $\delta_{m,q}=I_{n,m,q}'/q^m$ and use \thetag{58} and \thetag{63}.
\qed
\enddemo

\remark{Remark}
For $q=2$, one could use $I_{m,n,q}'$ with $m=(q^n-1)/(q-1)$
to derive a rationality criterion and an irrationality test for~$\gamma$
(but they are already contained in Theorems~2 and~3 as the case
$m=2^n-1$, since $I_{n,2^n-1,2}'=I_{n,2^n-1}$ and
$L_{n,2^n-1,2}'=L_{n,2^n-1}$). For $q>2$, this is not possible,
the required inequality $d_{q^n}I_{n,m,q}'<1$ not being valid,
since
$$
\lim_{n\to\infty}\frac{\log d_{q^n}}{q^n}=1
$$
(from the Prime Number Theorem) and the following observation
imply that $d_{q^n}I_{n,m,q}'$ tends to infinity with~$n$.
\endremark

\proclaim{Lemma 13}
If $q\ge3$, then
$$
\lim_{n\to\infty}\frac{\log I_{n,(q^n-1)/(q-1),q}'}{q^n}>-1.
\tag68
$$
\endproclaim

\demo{Proof}
By \thetag{64}, we have
$$
\split
\lim_{n\to\infty}\frac{\log I_{n,(q^n-1)/(q-1),q}'}{q^n}
&=\frac1{q-1}\lim_{n\to\infty}
\frac{\log I_{n,(q^n-1)/(q-1),q}'}{(q^n-1)/(q-1)}
\\
&>-\frac{1+\log(q+1)}{q-1}
\qquad\text{for}\quad q\ge2.
\endsplit
$$
Since the lower bound is increasing, and exceeds $-1$ when $q=4$,
this proves \thetag{68} for $q>3$. For $q=3$, we use \thetag{67}
to replace the bound by the exact value
$$
\frac12\log f_3(x_3)
=-0.90997390\hdots>-1,
$$
where $x_3=0.86304075\dots$ is the positive root of
$12-7x-8x^2$, which is the polynomial \thetag{62} with $q=3$.
\qed
\enddemo

\remark{Remark}
Alternatively, for $q$ sufficiently large \thetag{68} follows
from \thetag{67} and the formula
$$
\lim_{q\to\infty}f_q(x_q)=\frac1{2e},
$$
which in turn follows from the upper bound in \thetag{61}
and the fact that the limit, as $q\to\infty$, of the lower bound in
$$
\max_{0\le x\le1}f_q(x)
\ge f_q(x)\big|_{x=(1-1/q)^{1/q}}
=\biggl(q-\frac{1/q}{1-(1-1/q)^{1/q}}\biggr)
\cdot\biggl(1-\frac1q\biggr)^{q-1}
$$
is $1/(2e)$.
\endremark

\head
10. Concluding remarks
\endhead

Any polynomial $P_m(x)\in\Bbb Z[x]$ can be used
to generalize integral~\thetag{18} to base~$q$ as follows:
$$
\int_0^1\biggl(\frac q{1-x^q}-\frac1{1-x}\biggr)P_m(x)
\sum_{\nu=n+1}^\infty x^{q^\nu-1}\,\d x.
\tag69
$$
Integrals \thetag{37} and \thetag{40} are particular cases
of construction~\thetag{69}. The choice of $P_m(x)$,
the ``damping factor,'' is based on the requirements:
\roster
\item"(i)" to ensure the appearance of~$\gamma$,
and no other $\gamma_{j,q}$ from~\thetag{39},
in the corresponding linear forms;
\item"(ii)" to produce sufficiently good asymptotics (upper estimates).
\endroster
The choice $P_m(x)=(1-x)^m$ gives the best possible answer \thetag{38}
to feature~(ii); that is why we use it in
Sections~2--4 for $q=2$, and Section~8 for $q=3$.
The following lemma and the arguments in Section~7 show that
$(1-x)^m$ cannot be used for $q>3$, because of requirement~(i).

\proclaim{Lemma 14}
If $q>3$, the equality
$$
(1-x)^m=C+Q(x)\cdot\frac{1-x^q}{1-x}
\tag70
$$
is impossible for any $m>0$, constant~$C$ and polynomial $Q(x)$.
\endproclaim

\demo{Proof}
For $q>1$, equality \thetag{70} implies
the system of equalities $(1-\epsilon^j)^m=C$
for $j=1,2,\dots,q-1$, where $\epsilon$~is a primitive
$q$-th root of unity.
In particular, the value of $|1-\epsilon^j|$~is
the same for all $j=1,2,\dots,q-1$,
which is possible only when $q=2$ or $q=3$.
\qed
\enddemo

One can avoid this difficulty by introducing linear combinations of
powers of $1-x$. More precisely, there exist (integer) constants
$c_{m,0},c_{m,1},\dots,c_{m,q-1}$ such that
$$
c_{m,0}(1-x)^m+c_{m,1}(1-x)^{m-1}+\dots+c_{m,q-1}(1-x)^{m-q+1}
=C+Q(x)\cdot\frac{1-x^q}{1-x}.
$$
However, we cannot choose the left-hand side as the
polynomial $P_m(x)$, because the constants grow so much faster than~$m$
that requirement~(ii) fails.


\enddocument